\documentclass{article}
\usepackage[utf8]{inputenc}
\usepackage{amsmath}
\usepackage{amsfonts}
\usepackage{amssymb}
\usepackage{amsthm}
\usepackage{graphicx}
\usepackage{yhmath} 
\usepackage{mathrsfs}
\usepackage{subcaption}
\usepackage[ruled,linesnumbered]{algorithm2e}

\usepackage{xcolor}

\usepackage[left=2.50cm, right=2.50cm, top=3.00cm, bottom=3.00cm]{geometry}
\newtheorem{defi}{Definition}

\theoremstyle{definition}
\newtheorem{theorem}{Theorem}
\newtheorem*{theorem*}{Theorem}
\newtheorem{lemma}{Lemma}
\newtheorem{coro}{Corollary}
\newtheorem{prop}{Proposition}

\newcommand{\K}{\mathcal{K}}
\newcommand{\F}{\mathcal{F}}
\newcommand{\collapse}{\searrow}
\newcommand{\R}{\mathbb{R}}
\newcommand{\A}{\mathcal{A}}
\newcommand{\B}{\mathcal{B}}
\newcommand{\D}{\mathcal{D}}
\newcommand{\V}{\mathcal{V}}
\newcommand{\Q}{\mathcal{Q}}

\renewcommand{\H}{\mathcal{H}}
\renewcommand{\subset}{\subseteq}
\renewcommand{\L}{\mathcal{L}}
\renewcommand{\hat}[1]{\widehat{#1}}

\title{Survey:\\
{\Large  Optimal Discrete Morse Theory Simplification}}
\author{Francisco Martinez}
\date{December 2018}

\begin{document}
	
	\maketitle
	
	A central problem in topological data analysis and specifically in Computational Topology is that of computing the homology of a given simplicial complex. Said complexes can have arbitrary large number of simplices, as can happen, for example, if the space is the Rips-Vietoris or \v{C}ech complex of a large data cloud. As of the moment of this survey, the bottleneck in algorithms for computing homology comes from matrix multiplication, which has a complexity of $O(n^\omega)$, where $2<\omega<3$. Because of this reason, pre-processing the simplicial complex to get a smaller complex with the same homology groups, and then applying the homology algorithm to the smaller one, has been an active research topic in the last years. In this survey we discuss some recent papers that examine the complexity of this simplification via Discrete Morse Theory.\\
	
	We start with a short review of Discrete Morse Theory. Then we follow a paper by Joswig and Pfetsch \cite{JoswigPfetsch2006} to prove that finding an optimal Morse Matching is NP-hard. We then discuss a paper by Bauer and Rathod \cite{Bauer2018}, to show that even trying to find a near optimal approximation is also a hard problem. Finally we review a paper by Rathod, Bin Masood, and Natarajan \cite{Rathodetal2017}, describing a polynomial algorithm to find Morse Matching close to the optimal within a factor of $(D+1)/(D^2+D+1)$, where $D$ is the dimension of the simplicial complex. While this factor is far from 1, according to experiments done in \cite{Rathodetal2017}, the algorithm gets a better approximation ratio for the tested datasets.

	\section{Discrete Morse Theory}
	We start with a review of Discrete Morse Theory. Among the many equivalent ways to present this theory, we choose to do it using acyclic matches in Hasse Diagrams, since this way will be useful in the following sections. However, this is all equivalent to Forman's original formulation with Morse functions and discrete gradients (see \cite{Forman2002,Forman95}). \\
	
	\begin{defi}[Poset]
		A \textbf{partially ordered set} or \textbf{Poset}, is a pair $(P,\leq)$ where $P$ is a set and $\leq$ is an order relation on $P$, i.e. it satisfies:
		\begin{itemize}
			\item Reflexivity: $x\leq x$ for all $x\in P$.
			\item Transitivity: If $x\leq y$ and $y \leq z$ then $x\leq z$.
			\item Antisymmetry: If $x\leq y$ and $y\leq x$ then $x=y$. 
		\end{itemize}
	
	Given a Poset $(P,\leq)$, we denote by $\prec$ the \textbf{covering relation} in $P$, that is, $x\prec y$ if $x < y$ and there doesn't exist any $z\neq x,y$ such that $x\leq z\leq y$. 	
	\end{defi}
	
	\begin{defi}[Hasse Diagram]
		Given a Poset $(P,\leq)$, its Hasse Diagram (or graph) is the directed graph $\mathcal{H}$ whose vertices are all elements of $P$ and directed edges $(x,y)$ if and only if $y\prec x$. 
	\end{defi}
	
	\newpage
	\begin{defi}\leavevmode
		\begin{enumerate}
			\item Given a poset $(P,\leq)$, a \textbf{partial matching} in $P$, is a partial Matching in the underlying graph of its Hasse Diagram, i.e. it is a subset $M\subset P\times P$ such that $(a,b)\in M$ implies $a\prec b$ and each $a\in P$ belongs to at most one pair in $M$. 
			\item A partial matching on $P$ is called \textbf{acyclic} if there does not exist any cycle
			$$ a_1\prec b_1\succ a_2\prec b_2 \cdots a_n\prec b_n\succ a_1 $$
			with $n\geq 2$, $(a_i,b_i)\in M$ and all the $b_i$'s distinct.			
		\end{enumerate}
	\end{defi}

	Note that an equivalent condition for a matching $M$ to be acyclic, is that if we take the Hasse Diagram $\mathcal{H}$ of $P$ (with the directed edges going from a simplex to a face) and change the orientation of all the edges in $M$ (to point from a simplex to a coface), the new directed graph, which we will denote $\mathcal{H}_M$, has no directed cycles. 
	
	\begin{defi}
		If $\K$ is a simplicial complex (or more generally a CW complex), its face poset $\F(\K)$ is the poset whose elements are all non-empty faces of $\K$ with the partial order of inclusion $\subset$. In this case, the covering relation $\sigma\prec \tau$ means $\sigma$ is a co-dimension 1 face of $\tau$. We denote the Hasse diagram of $\F(\K)$ by $\H_\K$.
	\end{defi}
	
	\begin{defi}[Morse Matching]\leavevmode
		\begin{enumerate}
			\item If $\K$ is a simplicial complex (or more generally a CW complex), an acyclic partial matching on its face poset $\F(\K)$ (which as discussed above corresponds to an acyclic orientation of its Hasse diagram $\H_{\K}$) is a \textbf{Morse Matching} on $\K$. 
			\item Given a Morse Matching $M$ on $\K$, the unmatched simplices, i.e. $\K\setminus M$, are the \textbf{critical simplices} of $M$. 
		\end{enumerate}
	\end{defi}

	Suppose that $\K_2\subset \K_1$ are two simplicial complexes such that $\K_2$ has exactly two simplices $\sigma$ and $\tau$ not in $\K_1$, it must be then that $\sigma$ is a \textbf{free face} of $\tau$ (i.e. $\sigma$ is contained in the unique maximal face $\tau$). Then is easy to see that $\K_2$ is a deformation retract of $\K_1$ (as topological spaces). This combinatorial deformation is known as a \textbf{simplicial collapse}. 
	
	\begin{defi}\leavevmode
		\begin{enumerate}
			\item Let $\K_1$ and $\K_2$ be simplicial complexes. If we can get $\K_2$ from $\K_1$ by doing finitely many simplicial collapses, we say that $\K_1$ \textbf{collapses} to $\K_2$ and we denote it by $\K_1\collapse \K_2$. Note that in this case, $\K_1$ and $\K_2$ are homotopy equivalent.
			\item If $\K\collapse 0$ (i.e. $\K$ can be collapsed to a point) we say $\K$ is \textbf{collapsible.}
		\end{enumerate}
	\end{defi}
	
	Suppose now we know $\K_1\collapse \K_2$, so there exists a sequence of elementary collapses $(\sigma_1,\tau_1),(\sigma_2,\tau_2),\cdots, (\sigma_n,\tau_n)$ going from $\K_1$ to $\K_2$. We can consider these pairs as a matching on $\K_1$. Moreover, suppose $a_1\prec b_1\succ a_2\prec\cdots a_k\prec b_k\succ a_1$, $k\geq 2$ is a cycle in the matching, without loss of generality we may assume $(a_1,b_1)$ is the first pair in the cycle appearing in the collapsing order. But that would mean $b_k$ is still in the simplicial complex when we collapse $(a_1,b_1)$, and $a_1\prec b_k\neq b_1$. This is a contradiction, since $a_1$ needed to be free in order to be able to collapse it. Therefore, $M$ is an acyclic Morse Matching. We just proved the following statement:
	
	\begin{prop}\label{prop:collapsesMorse}
		Let $\K$ be a simplicial complex and $\mathcal{L}\subset\K$ a subcomplex. If $\K\collapse\mathcal{L}$, then $\K\setminus\mathcal{L}$ can be paired to get a Morse Matching, and such pairing is given by the elementary collapses. 
	\end{prop}
	
	The converse of this proposition is also true, more generally, even when the critical simplices do not form a simplicial complex, these can be used to determine the homotopy type of $\K$. This is the result of the Main Theorem of Discrete Morse Theory, which we state below without proof. For a complete presentation with proof, refer to \cite{Kozlov2008,Forman95}. 
	
	\begin{theorem}[Main Theorem of Discrete Morse Theory]\label{thm:Main-Thm-Morse}
		Let $\K$ be a simplicial complex, and let $M$ be a Morse Matching. Let $c_i$ denote the number of critical $i$-dimensional simplices of $\K$.
		\begin{enumerate}
			\item Let $\mathcal{L}\subset \K$ be a subcomplex such that $\K\setminus \mathcal{L}$ is a union of pairs in $M$. Then $\K\collapse\mathcal{L}$. In particular, if the critical simplices form a subcomplex $\K_c$, then $\K\collapse\K_c$.
			\item In general, $\K$ is homotopy equivalent to a CW complex with $c_i$ cells of dimension $i$. 
		\end{enumerate}
	\end{theorem}

	Discrete Morse Theory also provides a relationship between the Betti numbers (that is, the dimension of the Homology groups $H_\ast(\K,\mathbb{F})$, where $\mathbb{F}$ is a field) and the number of critical simplices of each dimension. This result is known as Morse Inequalities, which we here state again without proof. 
	
	\begin{theorem}[Morse Inequalities]\label{thm:MorseIneq}
		Let $\K$ be a simplicial complex and $M$ a Morse Matching on $\K$. Let $c_i$ denote the number of critical $i$-dimensional simplices of $\K$. Let $\mathbb{F}$ be a field, and denote by $\beta_i=\dim(H_i(\K,\mathbb{F}))$ the $i-$th Betti number. Then we have
		$$c_d-c_{d-1}+c_{d-2}-\cdots+ (-1)^dc_0\geq \beta_d-\beta_{d-1}+\beta_{d-2}-\cdots+(-1)^d\beta_0$$
		for each $d=0,1,2,\dots$. 
		
	\end{theorem}
	
	\paragraph*{Examples}
	\begin{enumerate}	
		\item Let $\K$ be the boundary of an $n$-dimensional simplex. So $\F(\K)$ consists of all non-empty proper subsets of $\{1,2,\dots,n+1\}$. Consider the matching $M$: $(S,S\cup\{1\})$ for all $S\subsetneq\{2,\dots,n+1\}$. Clearly this match is acyclic, so we get a Morse Matching. The only critical simplices are $\{1\}$ and $\{2,3,\dots,n+1\}$. The main theorem then implies that $\K$ is homotopy equivalent to a CW complex consisting of one vertex and one $n-1$-dimensional face. The only such CW complex is the sphere $S^{n-1}$, so $\K\simeq S^{n-1}$ as expected. 	
		
		\item Figure \ref{fig:morse-projectiveplane} shows a Morse Matching on a triangulation of the projective plane $\mathbb{R}P^2$. The pairs $(\sigma,\tau)$, with $\sigma\prec\tau$ are shown by arrows pointing from $\sigma$ to $\tau$. Note the vertex 1, the edge $e$ and the triangle $t$ are not paired, so these are the critical simplices. The main theorem then implies that $\mathbb{R}P^2$ is homotopy equivalent to a CW complex with one 0-cell, one 1-cell and one 2-cells. Sadly, there are many CW complexes with this number of cells. However, by studying more carefully the Morse Matching it is possible to recover the topology of $\mathbb{R}P^2$, see for example \cite{Forman2002}. We will not get into these details in this survey. 
		\begin{figure}
			\centering
			\includegraphics[width=0.3\linewidth]{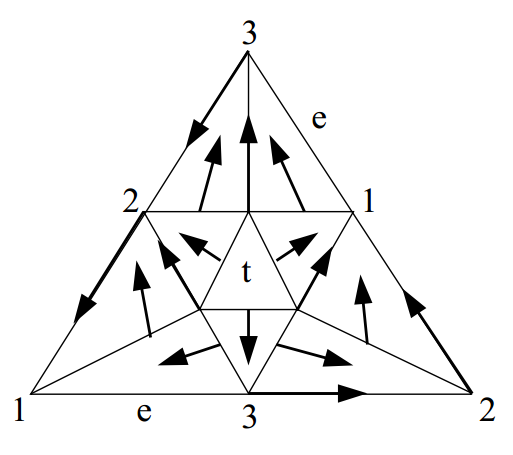}
			\caption{A Morse Matching on $\mathbb{R}P^2$. Taken from \cite{Forman2002}}
			\label{fig:morse-projectiveplane}
		\end{figure}
			
	\end{enumerate}

	These examples show how Theorem \ref{thm:Main-Thm-Morse} can be used to simplify a given simplicial complex, making it less expensive to compute, for example, homology or homotopy groups. It is thus a central problem of computational topology to find optimal (or near-optimal) Morse Matchings. Here, by optimal Morse Matching we mean the Morse Matching having the least critical simplices, or equivalently, having the most matched faces. 
	
	\section{Morse Matching is NP-hard}
	Because of its many applications, it is of vital importance to understand the complexity of computing optimal Morse Matchings of simplicial complexes, or at least for some special subcollections of complexes, like manifolds. Sadly, as Joswig and Pfetsch prove in \cite{JoswigPfetsch2006}, this problem is NP-hard. In this section we follow their proof of this fact.\\
	
	E\v{g}ecio\v{g}lu and Gonzalez proved in \cite{EgeciogluGonzalez1996} that the following problem on collapsibility of 2-dimensional complexes is NP-complete: Given a simplicial complex $\K$ embeddable in $\R^3$ and $k\in\mathbb{N}$, decide whether there exists a set $F$ of 2-faces with $|F|\leq k$, such that $\K\setminus F$ collapses to a 1-dimensional complex. By using proposition \ref{prop:collapsesMorse} and Theorem \ref{thm:Main-Thm-Morse}, we can restate this result as the following theorem.
	
	\begin{theorem}[E\v{g}ecio\v{g}lu and Gonzalez, 1996]\label{thm:collapsibility-problem}
		Given a connected 2-dimensional complex embeddable in $\R^3$, and $k\in\mathbb{N}$, it is NP-complete to decide whether there exist a Morse matching with at most $k$ 2-dimensional critical faces. 
	\end{theorem} 

	We now prove the equivalence between the statements.
	
	\begin{proof}[Proof of the equivalence]
		Suppose first that there is a set $F$ of 2-faces with $|F|\leq k$ and such that $\K\setminus F\collapse \mathcal{L}$, where $\mathcal{L}$ is a 1-dim. complex. Thus, by proposition \ref{prop:collapsesMorse}, $\K\setminus F \setminus\mathcal{L} = \K\setminus (F\cup\mathcal{L})$ can be paired in a Morse Matching. So the critical simplices are in $\mathcal{L}\cup F$, thus all the 2-dim. critical faces are in $F$. \\
		
		For the converse, suppose $M$ is a Morse Matching and $F$ is the set of 2-dim. critical faces, with $|F|\leq k$. Thus all 2-faces in $\K\setminus F$ appear in pairs $(\sigma_i,\tau_i)$ in the matching, where $\tau_i$ is 2-dimensional. Note here that at least one of the $\sigma_i$'s has to be free, otherwise we would be forced to have a cycle $\tau_1\succ\sigma_1\prec\tau_2\succ\cdots\prec\tau_r\succ\sigma_r\prec\tau_1$. Say thus $\sigma_1$ is free, so we can collapse $(\sigma_1,\tau_1)$. By an inductive argument, we can collapse all pairs $(\sigma_i,\tau_i)$ in some order, getting rid of all 2-dimensional faces in $\K\setminus F$. Thus $\K\setminus F$ collapses to a 1-dimensional complex as desired.  
	\end{proof}

	\begin{lemma}
		Let $\K$ be a connected simplicial complex. Let $\Gamma(M)$ be the graph obtained from the 1-skeleton of $\K$ by removing all edges matched with 2-faces in $M$ (note $\Gamma(M)$ contains all vertices). Then $\Gamma(M)$ is connected. 
	\end{lemma}

	\begin{proof}
		Assume $\dim(\K)\geq 2$, otherwise the result is immediate. Suppose $\Gamma(M)$ is disconnected. Let $N$ be one of its connected components and let $C$ be its cut edges (i.e. those edges in $\K$ with exactly one vertex in $N$). For each $\sigma\in C$, there exists a 2-face $\tau$ such that $(\sigma,\tau)\in M$. Say $\tau=abc$, $\sigma=ab$, $a\in N$ and $b\not\in N$, then either $ac$ or $bc$ will also be a cut edge in $C$, call it $\tau_2$. By continuing this process we produce a path $\sigma_1\prec\tau_1\succ\sigma_2\prec\tau_2\succ\sigma_3\prec\cdots$. Since there are finitely many simplices, this path must have a cycle, producing in this way a contradiction. 
	\end{proof} 
    
    The following result follows from the lemma we just proved by also realizing that on a connected 1-dimensional simplicial complex, we can construct a Morse Matching with unique critical vertex any desired vertex. This can be done via depth first search, which in turn can be computed in polynomial time. We skip the details, refer to \cite{JoswigPfetsch2006,Bauer2018}.
    \begin{lemma}\label{lemma:smaller-matching}
		Let $\K$ be a connected simplicial complex, let $p$ be a vertex of $\K$ and let $M$ be a Morse Matching on $\K$ with $c_0> 1$ critical vertices and $c$ critical simplices in total. Then there exists another Morse Matching $\tilde{M}$ on $\K$  with $p$ as the only critical vertex, $c-2(c_0-1)$ critical simplices in total, and the same number of $i-$dimensional critical simplices as $M$ for $i\geq2$. Moreover, this new matching can be computed from $M$ in polynomial time. 
	\end{lemma}

	We get now the desired result. For completion, we reproduce the proof given by Joswig and Pfetsch in \cite{JoswigPfetsch2006}.
	
	\begin{theorem}[Joswig and Pfetsch, 2006]
		Given a simplicial complex $\K$ and $c\in\mathbb{N}$, it is NP-complete to decide whether there exists a Morse matching with at most $c$ critical faces. Even if $\K$ is 2-dimensional and embeddable in $\R^3$. 
	\end{theorem}
	
	\begin{proof}
		Clearly this problem is NP (given a matching it can be checked in polynomial time if it is acyclic, and we can find the maximal one by exhaustion). Let $(\K,k)$ be an input for the collapsibility problem. We'll prove: there exists a Morse matching with at most $k$ critical 2-faces if and only if there exists a Morse matching with at most $g(k)=2(k+1)-\chi(\K)$ critical simplices in total (i.e. critical vertices, edges and faces). Here $\chi(\K)=\beta_0-\beta_1+\beta_2$ is the Euler Characteristic of $\K$, which can also be computed in terms of the number of faces in $\K$. Hence this problem is equivalent to the problem in Theorem \ref{thm:collapsibility-problem}, and so is NP-complete. \\
		
		Assume that $M$ is a Morse matching on $\K$ with at most $k$ critical 2-faces. Lemma \ref{lemma:smaller-matching} produces in polynomial time a new Morse matching $M'$, such that $c_0(M')=1, c_2(M')=c_2(M)$ and $c(M')\leq c(M)$. By the Main Theorem of Morse Theory \ref{thm:Main-Thm-Morse}, we have $\chi(\K)=c_0(M')-c_1(M')+c_2(M')$, from where we get 
		$$c(M')=2+2c_2(M')-\chi(\K)\leq 2(1+k)-\chi(\K).$$
		Conversely, suppose there is a Morse matching $M$ with at most $g(\K)$ critical simplices. Again lemma \ref{lemma:smaller-matching} produces $M'$, and so we get
		$$c_2(M')=c_2(M)=\frac{1}{2}\left(c(M)+\chi(\K)\right)-1\leq \frac{1}{2}\left(g(k)+\chi(\K)\right)-1=k$$
		 
	\end{proof}
	
	\section{Hardness of Approximating Optimal Morse Matching}
	In last section we discussed that finding an Optimal Morse Matching is an NP-complete problem. Since most of the time we don't need the optimal Morse Matching, but instead all we care about is to reduce the original simplicial complex, finding just a near optimal solution is good enough. It is then natural to ask if it is possible to find a polynomial time algorithm that produces near optimal Morse Matchings. On this section we discuss recent work by Bauer and Rasthod in \cite{Bauer2018} showing that even just looking for near-optimal solutions is already hard.\\
	
	When looking for the optimal Morse Matching it is equivalent to look for the one with most regular faces (paired simplices in the match) than for the one with least critical faces. Since we now are dealing with approximations, these turn out to be quite different problems. We now formally defined these problems. 
	
	\begin{defi}\leavevmode
		\begin{enumerate}
			\item \textbf{Max-Morse matching problem} (MaxMM): Given a simplicial complex $\K$, compute a Morse Matching that maximizes the number of matched simplices.
			\item \textbf{Min-Morse matching problem} (MinMM): Given a simplicial complex $\K$, compute a Morse Matching that minimizes the number of critical simplices. 
		\end{enumerate}
	\end{defi}

	\begin{defi}
		An \textbf{$\boldsymbol\alpha$-approximation} algorithm for an optimization problem is a polynomial time algorithm that for all instances of the problem produces a solution that is within an $\alpha$ factor of the optimal one. 
	\end{defi}
	
	Thus an $\alpha$-approximation algorithm for MaxMM, would be a polynomial time algorithm such that for each $\K$ simplicial complex, produces a Morse matching $M$ with $\dfrac{|M|}{|M_\text{opt}|}\geq\alpha$, where $M_\text{opt}$ is the optimal Morse Matching. Similarly, an $\alpha$-approximation algorithm for MinMM, would produce a Morse Matching $M'$ such that $\dfrac{c(M')}{c(M_\text{opt})}=\dfrac{n-2|M'|}{n-2|M_\text{opt}|}\leq \alpha$, here $n$ is the number of simplices in $\K$ and $c(M)$ are the critical simplices of $M$.

	\subsection{Hardness of MinMM}
	We now reproduce Bauer and Rasthod's proof (\cite{Bauer2018}) that unless P=NP, there is no $O(n^{1-\epsilon})$-approximation algorithm for MinMM, for any $\epsilon>0$, where $n$ is the number of simplices in the complex. 
	
	\begin{defi}[Amplified complex]
		Given a pointed simplicial complex $\K$ with $n$ simplices and an integer $c>0$, the \textbf{amplified complex} $\hat{\K}_c$ is the wedge sum of $n^{c-1}$ copies of $\K$.
	\end{defi}

	\begin{lemma}\label{lemma:amplified-complex}
		Given a simplicial complex $\K$ of size $n$ and $c\in\mathbb{N}^+$. Let $\hat{M}$ be a Morse Matching on the amplified complex $\hat{\K}_c$. Then
		\begin{enumerate}
			\item $\hat{\K}_c$ is collapsible if and only if $\K$ is collapsible. 
			\item If $\K$ is not collapsible, then $\hat{M}$ has more than $n^{c-1}$ critical simplices. 
		\end{enumerate}
	\end{lemma}

	\begin{proof}
		Suppose $\K$ is collapsible. By proposition \ref{prop:collapsesMorse}, there is a Morse Matching $M$ with only one critical vertex $q$. Let $p$ the distinguished point in $\K$ used to get $\hat{\K}_c$. By lemma \ref{lemma:smaller-matching}, we get the Matching $M_p$ on $\K$ with $p$ as its unique critical vertex, and we can copy this matching into all the copies of $\K$ to get a matching $\hat{M}_p$ of $\hat{\K}_c$. This is a Morse matching of $\hat{\K}_c$ with only one critical vertex, so by theorem \ref{thm:Main-Thm-Morse}, $\hat{\K}_c$ is collapsible.\\
		
		Suppose now $\hat{\K}_c$ is collapsible. So we can get a Morse Matching $\hat{M}$ with only one critical vertex. By using lemma \ref{lemma:smaller-matching}, we can transform this into $\hat{M}_p$, having as its unique critical vertex the distinguished point in $\hat{\K}_c$. We can then restrict this matching to a matching in one of the copies of $\K$, so $\K$ is also collapsible. \\
		
		For the second statement, suppose by way of contradiction that $\K$ is not collapsible and that $\hat{M}$ has $n^{c-1}$ or less critical simplices. By lemma \ref{lemma:smaller-matching}, we can suppose $\hat{M}$ has the distinguished point $p$ as its unique critical vertex, so the remaining critical simplices are of higher dimension. It must be then that one of the copies of $\K$ contains no other critical simplex than $p$. If we restrict $\hat{M}$ to this copy of $\K$, we would then get $\K$ to be collapsible, which is a contradiction. Therefore, $\hat{M}$ has less than $n^{c-1}$ critical simplices. 
	\end{proof}

	\begin{lemma}
		Let $0<\epsilon\leq1$. If there exists an $O(n^{1-\epsilon})$-approximation algorithm for MinMM ($n$ the number of simplices in the simplicial complex), then there exists a polynomial time algorithm for deciding whether a simplicial complex is collapsible or not. 
	\end{lemma}

	\begin{proof}
		Suppose such an approximation algorithm exists. In particular, there exist $p, N>0$ such that if $n\geq N$, the approximation ratio is at most $pn^{1-\epsilon}$. Let $c=\left\lfloor\frac{1}{\epsilon}+1\right\rfloor$. Let $\K$ be a connected simplicial complex with $n$ simplices, $n>\max\left(p^{1/(c\epsilon-1)},N\right)$, and construct the amplified complex $\hat{\K}_c$. The total number of simplices in $\hat{\K}_c$ is $\hat{n}=(n-1)n^{c-1}+1=n^c-n^{c-1}+1$. Consider now the following algorithm to decide the collapsibility of $\K$:
		Run the approximation algorithm to get a Morse Matching for $\hat{\K}_c$ (this can be done in polynomial time for $\hat{n}$, so is also polynomial time for $n$), and let $C_\text{approx}$ be the number of critical simplices of this matching. If $C_\text{approx}<n^{c-1}$ return $\K$ is collapsible, otherwise $\K$ is not collapsible. \\
		
		We now check this algorithm is correct. If $\K$ is not collapsible, by lemma \ref{lemma:amplified-complex} any Morse matching on $\hat{K}_c$ must have less than $n^{1-c}$ critical simplices, so the algorithm correctly says it is not. If $\K$ is collapsible, by the same lemma, $\hat{K}_c$ is also collapsible, and so its minimal Morse matching must have only 1 critical simplex, since the algorithm we used is a $p\hat{n}^{1-\epsilon}$ approximation, we get
		$$C_\text{approx}\leq p\hat{n}^{1-\epsilon}<p(n^c)^{1-\epsilon} < n^{c\epsilon-1}n^{c-c\epsilon}=n^{c-1}$$
		So the algorithm correctly says it is. 
	\end{proof}

	We now invoke a theorem by Tancer \cite{Tancer2016} regarding the collapsibility of 3 dimensional complexes. 
	
	\begin{theorem}[Tancer, 2016]
		It is NP-complete to decide whether a given 3 dimensional simplicial complex is collapsible. 
	\end{theorem}

	As a corollary we get the corresponding complexity for the Min-Morse matching problem. 
	
	\begin{theorem}
		For any $0<\epsilon\leq 1$, there doesn't exist a $O(n^{1-\epsilon})$-approximation algorithm for MinMM, unless P=NP. 
	\end{theorem}

	It is relevant to note here that while this means that a near optimal approximation to the general MinMM problem is NP-complete, even if we restrict to just 3 dimensional complexes, this gives no information about the restricted problem to 2 dimensional complexes. The case for just 2-dimensional complexes is, to the best of my knowledge, still an open question. 
	
	\subsection{Hardness of MaxMM}
	Not surprisingly, near optimal algorithms for MaxMM are also NP-hard, as Bauer and Rasthod also prove in \cite{Bauer2018}. To get this result, they construct an approximation preserving reduction from the 3OMAS problem to MaxMM, employing then a result by Newman \cite{Newman2001} that approximating 3OMAS is NP-complete yields the same complexity for the MaxMM problem. In this section we discuss the main ideas of this reduction, however we will skip the proofs of the more technical results.\\
	
	\begin{defi}[L-reduction]
		Let $\A$ be a maximization problem with objective function $m_\A$. For an instance $x$ of $\A$, $\normalfont\text{OPT}_\A(x)=\max_ym_A(x,y)$ is the optimal value of $m_\A$ for $x$, where the max is taken among all feasible solutions $y$ for $x$. Similarly, let $\B$ be another maximization problem with optimal value $\normalfont\text{OPT}_\B(x)$ for an instance $x$.
		
		An \textbf{L-reduction} from $\A$ to $\B$ with parameters $\mu$ and $\nu$ (positive constants), is a pair $(f,g)$ of polynomial time computable functions, satisfying:
		\begin{enumerate}
			\item For each instance $x$ of $\A$, $f(x)$ is an instance of $\B$.
			\item For each instance $x$ of $\A$, and each solution $y$ of $f(x)$, $g(x,y)$ is a solution of $x$. 
			\item $\normalfont\text{OPT}_\B(f(x))\leq \mu\text{OPT}_\A(x)$, for all instances $x$ of $\A$. 
			\item $\normalfont \text{OPT}_\A(x)-m_\A(x,g(x,y))\leq\nu\left[\text{OPT}_B(f(x))-m_\B(f(x),y)\right]$, for all instances $x$ of $\A$ and all solutions $y$ of $f(x)$. 
		\end{enumerate}
					
	\end{defi}
	
	The next theorem shows the usefulness of L-reductions to the problem in hand.
	
	\begin{theorem}[Vazirani \cite{Vazirani2003}]\label{thm:l-reductions}
		If there exists an L-reduction with parameters $\mu$ and $\nu$ from a maximization problem $\A$ to another maximization problem $\B$, and there is a $(1-\delta)$-approximation algorithm for $\B$, then there is a $(1-\mu\nu\delta)$-approximation algorithm for $\A$. 
	\end{theorem}

	\begin{defi}
		A directed graph $G$ is called an \textbf{oriented graph} if its underlying undirected graph is simple, i.e. it has no loops and no 2-cycles. 
		A \textbf{feedback arc set} is a subset of edges $E'$, such that $G'=(V,E\setminus E')\leq G$ has no directed cycles.  
		A \textbf{directed degree-3 graph} is a directed graph with total degree (adding in-degree and out-degree) at most 3. 
	\end{defi}

	\begin{defi}\leavevmode
		\begin{enumerate}
			\item \textbf{Maximal acyclic subgraph problem (MAS):} Given a directed graph $G$, find a maximal edge subset $E_\text{max}\subset E$ such that the graph $G_\text{max}=(V,E_\text{max})$ is a directed acyclic graph. 
			\item \textbf{Minimum cardinality Feedback arc set problem (minFAS):} Given $G=(V,E)$ a directed graph, find a minimal feedback arc set.
		\end{enumerate} 
	\end{defi}
	
	Clearly MAS and minFAS are complementary problems. We will also consider the MAS problem restricted to degree-3 graphs (3MAS), to oriented graphs (OMAS), and to oriented degree-3 graphs (3OMAS). We state the following result by Bauer and Rasthod without proof. 
	
	\begin{theorem}
		There exist an L-reduction with parameters $\mu=\nu=1$, from MAS to OMAS, and from 3MAS to 3OMAS. 
	\end{theorem}
	
	The following theorem by Newman \cite{Newman2001} gives the complexity for 3MAS. 
	
	\begin{theorem}[Newman,2001]
		It is NP-hard to approximate 3MAS within $\left(1-\frac{1}{126}\right)+\epsilon$, for any $\epsilon>0$. 
	\end{theorem}

	\begin{coro}\label{coro:complexity3OMAS}
		It is NP-hard to approximate 3MAS and 3OMAS to within $\left(1-\frac{1}{126}\right)+\epsilon$, for any $\epsilon>0$. 
	\end{coro}

	The main result Bauer and Rasthod prove in \cite{Bauer2018} is that there exists an L-reduction from 3OMAS to MaxMM with parameters $\mu=78$ and $\nu=\frac{1}{2}$. Applying then corollary \ref{coro:complexity3OMAS} and theorem \ref{thm:l-reductions}, produces the following theorem.
	
	\begin{theorem}[Bauer and Rasthod, 2018]\label{thm:maxmm_approx_hardness}
		It is NP-hard to approximate MaxMM within a factor of $\left(1-\frac{1}{4914}\right)+\epsilon$, for any $\epsilon>0$. 
	\end{theorem}

	In the next subsection we outline the main ideas of the L-reduction from 3OMAS to MaxMM.
	
	\subsubsection{L-Reduction form 3OMAS to MaxMM}
\begin{figure}
		\centering
		\begin{subfigure}[a]{0.45\textwidth}
			\centering
			\includegraphics[width=0.6\linewidth]{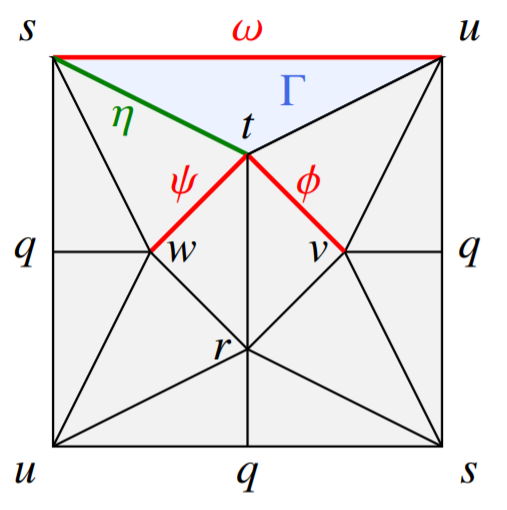}
			\caption{Modified Dunce hat $\D$, the vertices with the same name are identified. The distinguished simplices are $\Gamma, \omega,\psi,\phi,\eta$ and their vertices. }
			\label{fig:Dunce_modified}
		\end{subfigure}\hfill
		\begin{subfigure}[a]{0.45\textwidth}
			\centering
			\includegraphics[width=0.5\linewidth]{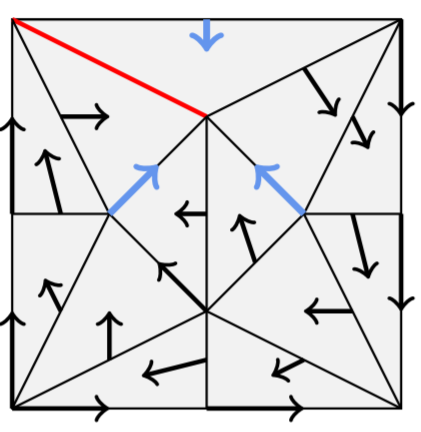}
			\caption{A Morse Matching on the Modified Dunce hat, such that only $s,t$ and $\eta$ are critical.}
			\label{fig:Dunce_MorseMatching}
		\end{subfigure}%
	\caption{Modified Duce hat $\D$, taken from \cite{Bauer2018}.}
	\end{figure}

	Given a degree-3 oriented graph $G$, to each subgraph $H$ we associate a simplicial complex $\K(G,H)$, and define $\K(G)=\K(G,G)$. To do such identification, for each edge $e$ in $E_H$ we consider a disjoint copy of the modified dunce hat $\D$ in figure \ref{fig:Dunce_modified} and denote it $\D_e$. We then identify some of the distinguished simplices of different copies of $\D$. While we won't discuss all the details of this identification, we do emphasize that only 0 and 1 dimensional simplices will get identified, and that this identification is done based on the in- and out-degrees of each of the vertices in $H$. Once the identification is done, the topological quotient space will be the desired $\K(H,G)$. Since the total degree of each vertex is at most 3, the identification rules can be carried in linear time. 
	
	\begin{prop}
		Given a degree-3 oriented graph $G$, and any subgraph $H$, $\K(H,G)$ is a simplicial complex that can be computed in polynomial time. Moreover, $\K(G,H)\leq\K(G)$ is a subcomplex.
	\end{prop}

	\begin{defi}
		We say a 2-dimensional simplicial complex $\K$ is \textbf{erasable} if it can be collapsed to a 1-dimensional complex. A maximal face $\tau$ of $\K$ is called \textbf{internal} if it contains no free faces. Denote by $\text{er}(\K)$ the minimal number of internal 2 faces that need to be removed to get an erasable complex.  
	\end{defi}
	
	\begin{defi}
		Given a 2-complex $\K$ and a subcomplex $\L\subset\K$, we say $\L$ is \textbf{erasable subcomplex of} $\K$ (through the Matching $M$) if there exists another subcomplex $\mathcal{C}\subset \K$ with $\K\collapse\mathcal{C}$ (induced by $M$) such that $\K\setminus\mathcal{C}\subset\L$ and $\K^{(2)}\setminus\mathcal{C}^{(2)}=\L^{(2)}$ (where $\K^{(2)}$ denotes the 2-simplices). 
	\end{defi}

	\begin{lemma}\label{lemma:erasability_properties}
		If $\L_1$ and $\L_2$ are erasable subcomplexes of a 2-complex $\K$, then so is $\L_1\cup\L_2$. 
	\end{lemma}

	We skip the proof of this lemma. We now establish some properties on $\K(H,G)$. Throughout the remaining of this section, $G=(V,E)$ will denote a degree-3 oriented graph, and $H=(V_H,E_H)$ a subgraph. 
	
	\begin{lemma}\label{lemma:lemma4.6}
		$H$ is acyclic if and only if $\K(G,H)$ is erasable.
	\end{lemma}
	
	We also skip the proof of this lemma.
	
	\begin{lemma}\label{lemma:lemma4.7}
		For any edge $e\in E$, $\D_e\setminus\{\Gamma_e\}$ is erasable in $\K(G)$.
	\end{lemma}

	\begin{proof}[Idea of the proof]
		Denote $M_e$ the Morse Matching depicted in figure \ref{fig:Dunce_MorseMatching} on $\D_e$. Then $M_e\setminus\{(\omega_e,\Gamma_e)\}$ erases $\D_e\setminus\{\Gamma_e\}$ in $\D_e$. Since the 1-simplices of $D_e$ paired with a 2-simplex in $M_e$ are not identified with those in any other $\D_{e'}$, the same matching shows $D_e\setminus\{\Gamma_e\}$ is erasable in $\K(G)$. 
	\end{proof}
	
	\begin{lemma}\label{lemma:lemma4.8}
		Let $C$ be a set of 2-simplices such that $\K(G)\setminus C$ is erasable, and let $F=\left\{f\in E:C\cap\D_f\neq\emptyset \right\}$. Then $F$ is a feedback arc set of $G$. 
	\end{lemma}

	\begin{proof}
		Note if $\sigma\in C$, then $\sigma\not\in\D_e$ for any $e\in E\setminus F$. Let $H=(V,E\setminus F)$ a subgraph of $G$, so $\K(G,H)\leq\K(G)\setminus C$ is erasable. A subcomplex of an erasable complex is also erasable, so lemma \ref{lemma:lemma4.6} implies $H$ is acyclic, so $F$ is a feedback arc set. 
	\end{proof}
	
	\begin{lemma}\label{lemma:ex-OPTminFAS}
		$\text{er}(\K(G))=\text{OPT}_\text{minFAS}(G)$
	\end{lemma}

	\begin{proof}[Proof omitting some details]
		Let $F\subset E$ be a minimal feedback arc set of $G$, and let $H=(V,E\setminus F)$, so $H$ is acyclic. 
		It can be shown that $\K(G)\setminus\{\Gamma_f|f\in F\}\collapse\K(G,H)$. By lemma \ref{lemma:lemma4.6}, $\K(G,H)$ is erasable, thus $\K(G)\setminus\{\Gamma_f|f\in F\}$ is erasable, and hence $\text{er}(\K)\leq|F|=\text{OPT}_\text{minFAS}(G)$.
		Suppose now $\text{er}(\K(G))<\text{OPT}_\text{minFAS}(G)$. Set $C$ be a minimal 2-simplices set such that $\K(G)\setminus C$ is erasable, so $\text{er}(\K(G))=|C|$. Let $F'=\{f\in E: C\cap\D_f\neq\emptyset\}$. By lemma \ref{lemma:lemma4.8}, $H=(V,E\setminus F')$ is acyclic, so $F'$ is a feedback arc set. Note that each 2-simplex in $C$ lies in at most one $\D_f$, so $|F'|\leq|C|$. Thus $|F'|\leq\text{er}(\K(G))<\text{OPT}_\text{minFAS}(G)$ but this is a contradiction. Hence the theorem follows.  
	\end{proof}

	Finally we state the next lemma showing the relationship between the topology of $\K(G)$ and $G$, however we will skip its proof. 
	
	\begin{lemma}
		$\K(G)$ is homotopy equivalent to $G$ (to the underlying graph).
	\end{lemma}
	
	The second function we need to have an L-reduction, needs to take instances of 3OMAS and solutions to the corresponding problem in MaxMM, to a solution to 3OMAS. We do this as follows: For a degree-3 oriented graph $G$, we compute a Morse Matching $M$ for $\K(G)$ and then define $F$ as follows:
	$$F=\left\{ f\in E|\ \ \exists\ \sigma\in\D_f,\ \sigma\text{ is a critical 2-simplex for }M \right\}$$
	By lemma \ref{lemma:lemma4.8}, this is a feedback arc set, thus the subgraph $H_M=(V, V\setminus F')$ is acyclic, and is the corresponding solution to 3OMAS on $G$. 
	
	\begin{lemma}\label{lemma:inequalities_objective_function}
		$m_\text{MaxMM}(\K(G),M)\leq n-2c_2(M)-\beta$ and $m_\text{3OMAS}(G,H_M)\geq |E|-c_2(M)$. Following the same notation as before, where $\beta=\beta_0+\beta_1+\beta_2$ is the sum of the Betti numbers of $\K(G)$, $M$ is a Morse matching and $c_i$ denotes the critical $i$-dimensional simplices of $M$. 
	\end{lemma}

	\begin{proof}
		Since $\K(G)\simeq G$, $\beta_2=0$ and $\beta_0=1$. By the Morse inequalities (theorem \ref{thm:MorseIneq}), $c_0\geq \beta_0$ and by some manipulations we get $$c=2c_2+\beta+2(c_0-\beta_0)\geq2c_2+\beta$$
		Thus $m_\text{MaxMM}(\K(G),M)\leq n-2c_2(M)-\beta$.\\
		
		In the construction of $H_M$, we remove at most as many edges as critical 2-simplices in $M$, so $|F'|\leq c_2$, so $m_\text{3OMAS}(G,H_M)\geq |E|-c_2(M)$.
	\end{proof}

	\begin{lemma}\label{lemma:Opt_maxmm}
		$\text{OPT}_\text{MaxMM}(\K(G))=n-2\text{OPT}_\text{minFAS}(G)-\beta$
	\end{lemma}

	\begin{proof}
		An optimal Morse Matching $M$ on $\K(G)$ will have $c_0=1$ and $c_2=\text{er}(\K(G))$. From lemma \ref{lemma:ex-OPTminFAS}, $\text{er}(\K(G))=\text{OPT}_\text{minFAS}(G)$. From the proof of lemma \ref{lemma:inequalities_objective_function}, when $c_0=\beta_0$ we have equality, thus: $$m_\text{MaxMM}(\K(G),M)=n-2c_2-\beta=n-2\text{minFAS}(G)-\beta.$$
	\end{proof}

	\begin{lemma}\label{lemma:parameter_mu }
		$\text{OPT}_\text{MaxMM}(\K(G))\leq 78\text{OPT}_\text{3OMAS}(G)$, so the L-reduction will have parameter $\mu=78$.
	\end{lemma}

	\begin{proof}
		In \cite{Vazirani2003}, they show that there always exists an acyclic subgraph with at least half of the edges of the original graph. So $\text{OPT}_\text{3OMAS}(G)\geq\frac{|E|}{2}$. The modified Dunce hat $\D$ has 39 simplices, so from the construction of $\K(G)$, $n\leq 39|E|$. Combining this two inequalities proves the result. 
	\end{proof}

	\begin{lemma}\label{lemma:parameter_nu}
		$\displaystyle\text{OPT}_\text{3OMAS}(G)-m_\text{3OMAS}(G,H_M)\leq\frac{1}{2}\Big[\text{OPT}_\text{MaxMM}(\K(G))-m_\text{MaxMM}(\K(G),M) \Big]$, so the L-reduction will have parameter $\nu=\frac{1}{2}$.
	\end{lemma}
	
	\begin{proof}
		This follows directly from combining lemmas \ref{lemma:inequalities_objective_function} and \ref{lemma:Opt_maxmm}.
	\end{proof}
	
	The construction $\K(G)$ from a oriented degree-3 graph $G$, and the solution $H_M$, an acyclic subgroup of $G$, obtained from a Morse Matching $M$ of $\K(G)$, constitute the L-reduction from 3OMAS to MaxMM with parameters $\mu=78$ and $\nu=1/2$ as lemmas \ref{lemma:parameter_mu } and \ref{lemma:parameter_nu} show. 
	
	\section{Approximation Algorithms for Max-Morse Matching}
	Theorem \ref{thm:maxmm_approx_hardness} establishes that it is NP-hard to find near Optimal solutions to the Max Morse Matching problem. However, if we content ourselves with approximate solutions within a smaller constant factor, it then becomes possible to do it in polynomial time. In this section, we discuss an algorithm developed by Rathod, Bin Masood and Natarajan in \cite{Rathodetal2017} in 2017, for which they are able to guarantee an approximation ratio of $\frac{D+1}{D^2+D+1}$, where $D$ is the dimension of the simplicial complex. While we won't discuss it here, we point out that they also produced an approximation algorithm for $D$-manifolds with an approximation ratio of $2/D$. \\
	
	Recall that given a simplicial complex $\K$, a Morse Matching $M$ is such that the oriented Hasse diagram $\H_M$ is acyclic. Where if $\sigma\prec\tau$ is in $M$ then the edge points from $\sigma$ to $\tau$, otherwise it points from $\tau$ to $\sigma$. Given a Hasse diagram $\H$, the \textbf{d-interface} of $\H$ is the subgraph consisting of all the nodes corresponding to dimension $d$ and $d-1$ simplices. Also note that all edges join simplices of co-dimension 1.
	
	\begin{defi}
		A simplex $\sigma$ of dimension $d$ is a \textbf{source node} for the $d$-interface, if it has only outgoing edges to $(d-1)$-simplices.\\
		A simplex $\tau$ of dimension $d-1$ is a \textbf{sink node} for the $d$-interface, if it has only incoming edges from $d$-simplices. 
	\end{defi}

	\begin{defi}
		An \textbf{up-edge} $\chi(\alpha,\beta)$ is an edge going from the $(d-1)$-simplex $\alpha$ to the $d$-simplex $\beta$, and is labeled $\chi$. We may sometimes ignore $\alpha$ and $\beta$ from the notation and only write $\chi$. A \textbf{down-edge}, goes from a $d$-dimensional simplex to a $(d-1)$-dimensional simplex. The \textbf{corresponding down-edge} to $\chi(\alpha,\beta)$ is the edge with orientation reversed, and is denoted by $\bar{\chi}(\beta,\alpha)$, or simply $\bar{\chi}$
	\end{defi}

	\begin{defi}
		If there is an edge $\chi_1(\alpha_1,\beta_1)$ followed by a down-edge $\bar{\chi}_2(\beta_1,\alpha_2)$ followed by an up-edge $\chi_3(\alpha_2,\beta_2)$, we say that $\chi_3$ is a \textbf{leading up-edge} of $\chi_1$.		
	\end{defi}

	\begin{defi}
		Given a $d$-dimensional simplex $\sigma$, the oriented edges between $\sigma$ and $(d-1)$-dimensional simplices, are the \textbf{facet-edges} of $\sigma$.
	\end{defi}

    \begin{lemma}
    	Let $M$ be a Matching on the Hasse graph of $\K$, and consider the oriented Hasse graph $\H_M$. Any directed cycle in $\H_M$ is completely contained in a $d$-interface, for some $1\leq d\leq\dim(\K)$.
    \end{lemma}

    \begin{proof}
    	Suppose there exists a cycle that passes through $\alpha$ and $\gamma$, where $\alpha$ and $\gamma$ are $(d-1)$-dimensional and $(d+1)$-dimensional, respectively. Thus, there must exist $\beta_1$ and $\beta_2$ $d$-dimensional such that $\chi_1(\alpha_1,\beta_1)$ and $\chi_2(\beta_2,\gamma)$ are up-edges. But to connect $\beta_1$ and $\beta_2$ at least one of them must be part of another up-edge, which is impossible since the up-edges come from a matching. This is then a contradiction, and the theorem follows. 
    \end{proof}

	\subsection{A ${}^{(D+1)}/{}_{(D^2+D+1)}$-factor Approximation Algorithm}
	
	The idea behind the algorithm is to start with a maximal Matching $M$ on the Hasse graph, the oriented Hasse graph $\H_M$ might have cycles, so the next step is to make sure we get an acyclic graph, by changing the orientation of some up-edges.\\
	
	The precise algorithm described by Rathod et al, is included in Algorithm \ref{algo:1}. Let us briefly review its complexity. The Hasse diagram can clearly be computed in polynomial time from the incidence information of the simplicial complex $\K$. Since all edges in $\H_\K$ only go from $d$-simplices to $(d-1)$ or $(d+1)$-simplices, $\H_\K$ is bipartite. A Maximum Matching on a bipartite graph can be found in polynomial time, for example by transforming the problem into a Maximum Flow problem and then applying Ford-Fulkerson algorithm (see \cite{Erciyes2018}, chapter 9.3). Clearly the \textbf{while} loop is bounded by the number of up-edges in $\H_M$, which, by definition, are the pairs in the matching $M$, hence the loop is executed at most $|M|\leq \frac{n}{2}$ times. We then just require to argue that the procedure \textbf{BFSComponent} is carried in polynomial time.\\
	
	The procedure \textbf{BFSComponent} is included in Algorithm \ref{algo:bfs}. Note that the \textbf{while} loop is executed as long as there is some edge in the queue $\Q$, and note that once an edge is dequeued, it never enters the queue again. Thus $\Q$ will at most include all the up-edges of $\H_M$, which, as we argued before, are at most $n/2$. Hence, the \textbf{while} loop is executed at most $n/2$ times. The inner \textbf{foreach} loop is also executed only for up-edges, so once again is bounded by $n/2$. Finally, at each stage we need to check if the graph with edges $C\cup$\textbf{facetEdges}$(\beta_i)$ has cycles, this can be done with depth first search since any loop created would need to include the new up-edge, so it can be done in polynomial time (see \cite{Erciyes2018}, chapter 6.3).
	
	\begin{algorithm}
		\SetKwFunction{BFSComponent}{BFSComponent}
		\KwData{Simplicial complex $\K$}
		\KwResult{$\H_\V$, an acyclic matching on the Hasse graph, based on a Morse Matching $\V$ on $\K$.}
		\vspace{1em}
		Construct the Hasse graph $\H_\K$ of $\K$\;
		Perform a Max-Cardinality Matching $M$ on $\H_\K$\;
		Let $\H_M$ be the Hasse graph with orientation based on the Matching $M$\;
		Initialize the graph $\H_\V$ with all the vertices of $\H_\K$ and edges $E(\H_\V)\leftarrow\emptyset$.

		\While{$\exists\chi\in E(\H_M)$ an up-edge}{
			$C_\chi\leftarrow$  \BFSComponent{$\H_M,\chi$}\;
			$E(\H_M)\leftarrow E(\H_M)\setminus C_\chi$\;
			$E(\H_\V)\leftarrow E(\H_\V)\cup C_\chi$\;
		}
        $E(\H_\V)\leftarrow E(\H_V)\cup E(\H_M)$\;
        \Return $\H_V$\;
        \caption{Frontier Edges Algorithm}\label{algo:1}
	\end{algorithm}
	
	\begin{algorithm}
		 \SetKwProg{proc}{Procedure}{}{}
		\SetKwFunction{BFSComponent}{BFSComponent}
		\proc{\BFSComponent{$\H_M,\chi$}}{
		      $C\leftarrow\emptyset$\;
		      Initialize a queue $\Q$\;
		      \textbf{enqueue}($\Q,\chi$)\;
		      \While{$\Q$ is not empty}{
		         $\chi_0(\alpha_0,\beta_0)\leftarrow$\textbf{deque}($\Q$)\;
		         $C\leftarrow C\cup\text{\textbf{facetEdges}}(\beta_0)$\;
		         \ForEach{leading up-edge $\chi_i(\alpha_i,\beta_i)$ of $\chi_0$}{
		                 \eIf{the graph with edges $C\cup\text{\textbf{facetEdges}}(\beta_i)$ has cycles}{
		                 reverse $\chi_i$ in $\H_M$\;
		                 $C\leftarrow C\cup\text{\textbf{facetEdges}}(\beta_0)$\;
	                 }{\textbf{enqueue}($\Q,\chi_i$)}}
	             }
		      \Return $C$
	      }
      \caption{BFS Component}\label{algo:bfs}
	\end{algorithm}
	
    \begin{lemma}
    	The graph with edges in an edge-component $C=$\textbf{BFSComponent}($\H_M,\chi_0$), is acyclic.    	
    \end{lemma}

    \begin{proof}
    	Note first that $C$ is contained in a $d$-interface, for some $d$. This since all leading up-edges $\chi_i(\alpha_i,\beta_i)$ to an edge $\chi$, are in the same $d$-interface, and so are all the edges in \textbf{facetEdges}($\beta_i$). Also note that all cycles in a $d$-interface are of the form $\alpha_0\beta_0\alpha_1\beta_1\cdots\alpha_k\beta_k\alpha_0$, where $\chi_{i+1}(\alpha_{i+1},\beta_{i+1})$ is a leading up-edge to $\chi_i(\alpha_i,\beta_i)$. In the \textbf{BFSComponent} procedure, when examining the up-edge $\chi_i(\alpha_i,\beta_i)$, all the other edges to which it leads to were already examined, so if it is checked that including \textbf{facetEdges}($\beta_i$), doesn't create any cycle, that is still the case when we dequeue $\chi_i$, so it can be added as an up-edge. If on the contrary, the check showed there was a cycle, we reverse $\chi_i$ and add \textbf{facetEdges}$(\beta_i)$ to $C$. If we do this, $\beta_i$ is not part of any up-edge anymore, and thus it can't be part of any cycle, similarly $\alpha_i$ is not part of any up-edge, so can't be in a cycle. So adding \textbf{facetEdges}$(\beta_i)$ produces no cycle in this case either. Therefore, $C$ has no cycles.     	
    \end{proof}

    \begin{lemma}
    	The output graph $\H_\V$ is acyclic.
    \end{lemma}

	\begin{proof}
		Following the proof in \cite{Rathodetal2017}, we prove this fact by induction on the sequential addition of the edge-components.
		\begin{itemize}
			\item \textbf{Base Step:} When including just one edge-component $C_1$, we get an acyclic graph, since $C_1$ is acyclic.
			
			\item \textbf{Inductive Step:} Suppose that the graph remains acyclic after the addition of the edges in the edge-component $C_i$. We need to prove it is still acyclic when we add the edge-component $C_{i+1}$. Suppose there is a cycle, so it is inside a $d$-interface, and must include at least one up-edge in $C_{i+1}$. Then it must also exists a down edge $(\beta_0,\alpha_1)$ in $C_{j_k}$, for some $j_k\leq i$, such that $\chi_1(\alpha_1,\beta_1)$ is an up-edge in $C_{i+1}$. But that would mean that $\alpha_1$ was reachable while computing the edge-component $C_{j_k}$ and thus $\chi_1$ was already included in $C_{j_k}$. Since the edge-components are disjoint, this is a contradiction. Therefore, there are no cycles after adding $C_{i+1}$.
		\end{itemize}
		When we reach line 9 of the algorithm \ref{algo:1}, $\H_\V$ is still acyclic. Finally, at line 10 we add all the remaining edges, in the facet-edges of simplices which were either unmatched by $M$ or were matched to their cofacets. In both cases, they are source nodes within the $d$-interface, and thus are not part of any cycle. Therefore, $\H_\V$ is acyclic.   
	\end{proof}

	\begin{theorem}\label{thm:algo_gives_morse}
		The output graph $\H_\V$ of Algorithm \ref{algo:1} is based on a Morse Matching on $\H_\K$.
	\end{theorem}

    \begin{proof}
    	Last lemma shows $\H_\V$ is acyclic, and all its up-edges come from the original matching $M$ computed on $\H_M$. Thus, all that need to be shown is that $\H_\V$ indeed includes all edges in the Hasse graph of $\K$, $\H_{\K}$ with some orientation. For this, we prove that for all $d$-simplices $\beta$, its faced-edges are included in $\H_\V$ at some time. We have two cases:
    	\begin{enumerate}
    		\item Case 1: If $\beta$ is matched to one of its facets in $M$. The up-edge incident to $\beta$ is examined at some time in the procedure \textbf{BFSComponent}, and at that moment all the facets of $\beta$ are included in $\H_\V$ (either keeping the up-edge or reversing it).
    		\item Case 2: If $\beta$ is unmatched or was matched with one of its cofacets, none of it edge-facets was reached by leading up edges to an up-edge, so they are all added in line 10 of Algorithm \ref{algo:1}. 
    	\end{enumerate}
    \end{proof}
	
	We now prove the approximation ratio of this algorithm. First we establish some notation. Recall that at a given step inside the \textbf{BFSComponent} procedure, we look at an up-edge and decide whether to include it in that direction or reversed. Once we've done this, we say that edge is \textbf{classified}. If it was decided to leave it as an up-edge, we call it \textbf{forward}, and if it is reversed we call it \textbf{backward}. Let the iterator $i$ be the number of up-edges already classified, let $F_i$ be the number of forward edges and $B_i$ be the number of backward edges, after classifying $i$ up-edges. If $\chi$ is an up-edge leading to one of the classified up-edges, that is not yet classified, we call it a \textbf{frontier edge}. Denote by $Z_i$ the number of frontier edges at step $i$. Note that at each step, the ratio $\dfrac{F_i}{F_i+B_i+Z_i}$ is the worst ratio of forward edges that we can have on the edges discovered so far (when all the frontier edges end up being classified as backward). 

	\begin{lemma}
		Let $C$ be an edge-component contained in the $d$-interface. If $U$ is the number of up-edges in $C$ as originally matched in $M$, and $F$ are those that are classified by the procedure \ref{algo:bfs} as forward, then $$\frac{F}{U}\geq\frac{d+1}{d^2+d+1}$$
	\end{lemma}
	
	\begin{proof}
		We proceed by induction on $i$ (the number of edges classified so far).
		\begin{itemize}
			\item \textbf{Base Case:} Let $\chi_0$ be the first up-edge visited by the procedure. Note that any cycle on $\H_M$ has minimum length 6: if $\sigma_1, \sigma_2\prec\tau$ are two different facets of $\tau$, by a cardinality argument, $\tau=\sigma_1\cup\sigma_2$ as sets, so they can't both be facets of a second $\tau_2$. Because of this, cycles can only appear after two iterations, so all edges will be classified as forward at $i=1,2$. Suppose $\chi_0$ has $K$ leading up-edges, and that each of them has $j_k$ leading up-edges, for $1\leq k\leq K$. Thus, after the second step we will have: $F_2=K+1$, $B_2=0$ and $Z_2=\sum_{k=1}^{K}j_k$.
			The worst ratio thus occurs when $K=j_k=d$ for all $k$, so
			$$\frac{F_2}{F_2+B_2+Z_2}\geq\frac{d+1}{d^2+d+1}$$
			
			\item \textbf{Induction Step:} Supposing that at time $i$, we have the ratio $\dfrac{F_i}{F_i+B_i+Z_i}\geq\dfrac{d+1}{d^2+d+1}$, we will prove this inequality still holds for step $i+1$. At time $i+1$, if the next up-edge is classified as backwards, we will have $F_{i+1}=F_i$, $B_{i+1}=B_i+1$ and $Z_{i+1}=Z_i-1$, since no more frontier edges are added. In this case, the ratio remains the same. If instead the new edge is classified as forward, $F_{i+1}=F_i+1$, $B_{i+1}=B_i$ and $Z_{i+1}\leq Z_1+d-1$, since at most the new forward edge will have $d$ leading up-edges. Since $1/d\geq (d+1)/(d^2+d+1)$ we have
			$$ \frac{F_{i+1}}{F_{i+1}+B_{i+1}+Z_{i+1}}\geq\frac{F_i+1}{(F_i+1)+B_i+(Z_i+d-1)}\geq\frac{d+1}{d^2+d+1}$$
		\end{itemize}
	After the last step $k$ of the procedure, we classify all the up-edges in that component of $\H_M$, so $U=F_k+B_k+Z_k$ and $F=F_k$, proving the inequality. 
	\end{proof}

	\begin{theorem}
		Algorithm \ref{algo:1} computes a $\frac{D+1}{D^2+D+1}$-approximation for the Max Morse Matching Problem, where $D$ is the dimension of the simplicial complex. 
	\end{theorem}
	
	\begin{proof}
		From Theorem \ref{thm:algo_gives_morse}, the Hasse diagram $\H_\V$ is based on a Morse Matching. The obtained Morse Matching is a subset of the original maximal matching $M$, made of the up-edges still in $\H_\V$. By the previous lemma, each execution of the \textbf{BFSComponent} procedure preserves a $\frac{d+1}{d^2+d+1}\geq\frac{D+1}{D^2+D+1}$ fraction of the up-edges, so at the end of algorithm \ref{algo:1}, we have $\dfrac{D+1}{D^2+D+1}|M|$ up-edges. Clearly $\text{OPT}_\text{MaxMM}(\K)\leq 2|M|$, so we have a ratio of
		$$\frac{2|\V|}{\text{OPT}_\text{MaxMM}(\K)}\geq \frac{2|\V|}{2|M|}\geq\frac{D+1}{D^2+D+1}$$
	\end{proof}
	
	\subsection{Other Algorithms}
	By exploiting the geometry of 2-simplicial complexes, Rahod, Bin Masood, and Natarajan, managed to prove a better approximation factor of Algorithm \ref{algo:1} for that case. Thus, Algorithm \ref{algo:1} provides a $5/11$ approximation ratio for 2-dimensional simplicial complexes. Using similar ideas, they also describe a better algorithm to deal with general $D$-manifolds, $D\geq 3$, for which they get an approximation ratio of $2/D$. 
	
	\subsection{Experimental Results}
	In the experiments run by Rathod, Bin Masood, and Natarajan, their algorithms actually produce better results that the guaranteed by the $(D+1)/(D^2+D+1)$ and $2/D$ ratios. They also compare their algorithms with some Heuristic algorithms based on reductions and coreductions, introduced in \cite{Harker2010}. We describe this heuristics here for completion, however we do not prove that they work.\\
	
	Suppose we are given a simplicial complex $\K$. 
	\begin{itemize}
		\item \textbf{Coreduction}. Perform the following steps on $\K$ until it is empty:
		\begin{enumerate}
			\item if there exists a simplex $\alpha$ with a free coface $\beta$, include the pair $(\alpha,\beta)$ to the Morse Matching, and delete them from $\K$. 
			\item if there doesn't exist any simplex with a free coface, select the lowest dimensional simplex $\gamma$. Make $\gamma$ critical and delete it from $\K$.
		\end{enumerate}
	\item \textbf{Reduction}. Perform the following steps on $\K$ until it is empty:
	\begin{enumerate}
		\item if there exists a simplex $\beta$ with a free face $\alpha$, include the pair $(\alpha,\beta)$ to the Morse Matching, and delete them from $\K$.
		\item if there doesn't exist any simplex with a free face, select the highest dimensional simplex $\gamma$. Make $\gamma$ critical and delete it from $\K$. 
	\end{enumerate}
	\end{itemize}
	
	Rathod, Bin Masood, and Natarajan's experiments show that the Coreduction heuristic is usually the best, however there doesn't exist any theoretical guarantee of its efficiency. On the other hand, the $2/D$ algorithm provides almost as good results, with the advantage of having theoretical bounds. Please refer to \cite{Rathodetal2017} for the details. 
	
	\section{Conclusions}
	
	Discrete Morse Theory is a very versatile and useful tool in Computational Topology and Topological Data Analysis because of its power to reduce the size of simplicial (and in general CW) complexes. Experiments show that Heuristics like Coreduction and algorithms like Algorithm \ref{algo:1} and the others described in \cite{Rathodetal2017} usually do a good job on the test datasets. However, the best approximation ratios so far for the Maximum Morse Matching problem are $(D+1)/(D^2+D+1)$ for $D$-simplicial complexes and $2/D$ for $D$-manifolds ($D\geq3$) as described in the last section.\\
	
	Finding the optimal Morse Matching is an NP-complete problem, as discussed in Section 2, following the work of Joswig and Pfetsch in \cite{JoswigPfetsch2006}. And even approximating either Max Morse Matching or Min Morse Matching with a near optimal ratio is NP-hard, as Bauer and Rathod show in \cite{Bauer2018}, and as we discuss in Section 3.

	\bibliographystyle{acm}
	\bibliography{References}
	
\end{document}